\documentclass[12pt,reqno]{amsart}
\usepackage{amsthm,amsmath,amssymb,amsfonts,amscd,verbatim,latexsym}

\newtheorem{Theorem}{Theorem}[section]
\newtheorem{Lemma}[Theorem]{Lemma}
\newtheorem{Proposition}[Theorem]{Proposition}
\newtheorem{Corollary}[Theorem]{Corollary}
\newtheorem{Remark}[Theorem]{Remark}

\newtheorem{Example}[Theorem]{Example}
\newtheorem{Problem}[Theorem]{Problem}

\newcommand{\complex}{\mathbf C}   
\renewcommand{\P}{\mathbb P}       
\renewcommand{\O}{\mathcal O}      
\newcommand{\R}{\mathcal R}        
\newcommand{\U}{\mathcal U}        
\newcommand{\G}{\mathcal G}        
\newcommand{\F}{\mathcal F}        
\newcommand{\J}{\mathbf J}         
\newcommand{\I}{\mathbb I}         
\newcommand{\Ev}{\mathcal E}       
\newcommand{\Res}{\text{Res}} 
\newcommand{\W}{\mathbb W}         
\newcommand{\ord}{\text{ord} \;}   
\newcommand{\sympol}{\mathfrak p}  
\newcommand{\ux}{\mathbf x}        
\newcommand{\uy}{\mathbf y}        
\newcommand{\polar}{\uy \, \partial_{\ux}} 
\newcommand{\Sym}{\text{Sym}}      

\newcommand{\cons}{{\mathbf k}\,}   
\newcommand{\ra}{\rightarrow}
\newcommand{\lra}{\longrightarrow}
\renewcommand{\phi}{\varphi}
\newcommand{\Gordan}[9]{ 
\left(\begin{array}{ccc} #1 & #2 & #3 \\ #4 & #5 & #6 \\ #7 & #8 & #9 
\end{array} \right)}             
\newcommand{\br}[2]{(#1 \, #2)}  
\newcommand{\po}[2]{#1^{\langle #2 \rangle}} 
\newcommand{\demo}{\noindent {\sc Proof.}\;}
\newcommand{\Hom}{\text{Hom}}

\begin{document}
\title{On the invariant theory of the B{\'e}zoutiant} 
\author{Jaydeep V.~Chipalkatti}
\maketitle 

\parbox{12cm}{ \small 
{\sc Abstract.} 
We study the classical invariant theory of the B{\'e}zoutiant $\R(A,B)$ for a 
pair of binary forms $A,B$. It is shown that $\R(A,B)$ is determined 
by the first two odd transvectants $M = (A,B)_1$, $N = (A,B)_3$, and 
one can characterize the forms $M,N$ which can arise from some $A,B$. 
We give formulae which express the higher odd transvectants 
$(A,B)_5,(A,B)_7$ in terms of $M$ and $N$. We also describe 
a `reduction formula' which recovers $B$ from $\R(A,B)$ and $A$.}

\bigskip 

\parbox{12cm}{\small
} 

\bigskip 

\section{Introduction} 
We begin by recalling the construction of the B{\'e}zoutiant of 
two binary forms. Let 
$\ux = (x_0,x_1), \uy=(y_0,y_1)$ be pairs of variables, and 
write $\omega = x_0 y_1 - x_1y_0$. If $A,B$ are (homogeneous) forms 
of order $d$ in $\ux$, their B{\'e}zoutiant is defined to be 
\[ \R(A,B) = \frac{A(x_0,x_1)B(y_0,y_1) - B(x_0,x_1)A(y_0,y_1)}{\omega}. \] 
Since $\R$ is symmetric in $\ux$ and $\uy$ of order $d-1$ in each,
it can be seen as a quadratic form over the vector space of 
order $d-1$ binary forms. 

If $V = \text{Span} \{x_0,x_1\}$, then the construction of $\R$ corresponds 
to the isomorphism of $SL(V)$-representations 
\[ \wedge^2 \, \Sym^d \, V \lra 
 \Sym^2(\Sym^{d-1} \, V), \quad 
A \wedge B \lra \R(A,B). \] 
It is easy to see that 
\[ \R(\alpha \, A + \beta \, B,\gamma \, A + \delta \, B) = 
(\alpha \, \delta - \beta \, \gamma) \, \R(A,B), \] 
i.e., up to a scalar, $\R$ depends only on the 
pencil spanned by $A,B$ (denoted $\Pi_{A,B}$). Conversely, 
$\R$ determines the pair $(A,B)$ up to a unimodular transformation. 

B{\'e}zoutiants have been principally studied 
for their use in elimination theory 
(e.g., see \cite{Householder1} or \cite[vol.~I, \S 136 ff]{Netto}). 
In contradistinction, our interest lies in their invariant theoretic 
properties (understood in the sense of Grace and Young \cite{GrYo}). 

\subsection{A summary of results} 
In section \ref{section.prelim} we recall some fundamental 
facts about transvectants. 
We will show that $\R(A,B)$ admits a `Taylor series' in 
$\omega$ as follows: 
\[ \R(A,B) = c_0 \, T_1^\sympol + c_1 \, \omega^2 \, T_3^\sympol + 
   c_2 \, \omega^4 \, T_5^\sympol + \dots, \] 
where 
\begin{itemize} 
\item $T_{2r+1}$ denotes the $(2r+1)$-th transvectant of 
$A,B$, 
\item 
$\sympol$ denotes the operation of symmetric polarization, and 
\item 
$c_r$ are rational constants dependent on $d$ and $r$. 
\end{itemize} 

{\sl Hence, from our viewpoint, a study of $\R(A,B)$ will be 
tantamount to a study of the odd 
transvectants $\{T_{2r+1}:r \ge 0\}$ of $A$ and $B$.}

\smallskip 

In section \ref{section.wr.ode}, we formulate a second order 
differential equation derived from $T_1,T_3$ whose solution 
space is $\Pi_{A,B}$. This shows that the terms of degree $\le 2$ in 
the Taylor series implicitly determine those of higher degree. The former 
cannot be chosen arbitrarily, and we give an algebraic 
characterization of terms which can so appear. Specifically, we construct a 
set of joint covariants $\Phi_0,\dots,\Phi_d$ with arguments $M,N$, 
with the following property: 

{\sl There exist $A,B$ such that 
$M = (A,B)_1, N = (A,B)_3$, if and only if 
$\Phi_0(M,N) = \dots = \Phi_d(M,N) = 0$.}

\smallskip 

We have remarked earlier that $\R$ determines $\Pi_{A,B}$.  
Hence, given $A$ and $\R$, the form $B$ is determined up to an 
additive multiple of $A$. In section \ref{section.gen.red}, we 
give an equivariant formula for $B$ in terms of $A$ and $\R$. 
This is called a `generic reduction formula', 
in analogy with a device introduced by D'Alembert in the theory of 
differential equations. 

In section \ref{section.t5t7}, we use the classical Pl{\"u}cker relations 
to describe formulae which calculate $T_5,T_7$ from a 
knowledge of $T_1$ and $T_3$. The question of a formula in 
the general case of $T_{2r+1},r \ge 4$ is left open. 
Three more open problems (with some supporting examples) are 
given in section~\ref{section.open}. 

\medskip 

{\sc Acknowledgements.} I thank my colleage 
A.~Abdesselam for several helpful conversations, and the University 
of British Columbia for financial support. 

\section{Preliminaries} \label{section.prelim}
We will heavily use \cite{GrYo} as a standard reference for classical 
invariant theory. Glenn's treatise \cite{Glenn} covers substantially 
the same ground. In particular, we assume some familiarity with 
transvectants, covariants, and the symbolic calculus. 
A more recent exposition of this material is given in \cite{Olver}. 
Basic facts about the representation theory of $SL_2$ can be found 
in \cite{FH,Springer1,Sturmfels}. 

The base field is throughout $\complex$. 
A form will always mean a homogeneous polynomial in $\ux$. 
By contrast, an $\ux\uy$-form will involve both sets of variables, 
and will be homogeneous in each set. The $\ux$-degree of a form 
will be called its order (to avoid conflict with \cite{GrYo}). 
The order of an $\ux\uy$-form is a pair of integers. 

The letter $\cons$ will stand for an 
unspecified nonzero constant. 

\subsection{$SL_2$-modules} 
Let $V$ be a $\complex$-vector space of dimension two with the 
natural action of $SL(V)$. We write $S_e$ for the 
symmetric power representation $\text{Sym}^e \, V$, and 
$S_e(S_f)$ for $\text{Sym}^e(\text{Sym}^f \, V)$ etc. 

The $\{S_e: e \ge 0\}$ are a complete set of finite dimensional 
irreducible $SL(V)$-modules. By complete reducibility, 
each finite dimensional $SL(V)$-module is isomorphic to a 
direct sum of the $S_e$. 
If $\{x_0,x_1\}$ is a basis of $V$, then an element of $S_e$ 
is a form of order $e$ in $\ux$. We identify 
the projective space $\P^e$ with $\P S_e$, and write 
$A \in \P^e$ for the point represented by a (nonzero) form $A$. 
By convention, $S_e =0$ if $e < 0$. 

\subsection{Transvectants} 
For integers $e,f \ge 0$, we have a decomposition of $SL(V)$-modules 
\begin{equation} 
S_e \otimes S_f \simeq \bigoplus\limits_{r=0}^{\min\{e,f\}} 
S_{e+f-2r}. 
\label{clebschgordan} \end{equation}
If $E,F$ are forms of orders $e,f$, the image of the projection 
of $E \otimes F$ in the $r$-th summand is called their $r$-th transvectant, 
denoted $(E,F)_r$. It is a form of order $e+f-2r$, whose coefficients are 
linear in the coefficients of $E$ and $F$. 
In coordinates, it is given by the formula 
\begin{equation}
(E,F)_r = \frac{(e-r)!(f-r)!}{e!f!}  
\sum\limits_{i=0}^r \, (-1)^i \binom{r}{i} \, 
\frac{\partial^{\,r} E}{\partial x_0^{r-i} \, \partial x_1^i } \; 
\frac{\partial^{\,r} F}{\partial x_0^i \, \partial x_1^{r-i}} 
\label{trans.defn} \end{equation}
(The initial scaling factor is conventional, some authors choose it 
differently.) In particular $(E,F)_0 = E \, F$, and 
$(E,F)_1 = \cons \times$ Jacobian$(E,F)$. Note that 
\begin{eqnarray} 
(E,F)_r & = & (-1)^r(F,E)_r, \label{trans.1} \\
(E,F)_r & = & 0 \qquad \text{for $r > \min\{e,f\}$}. \label{trans.2} 
\end{eqnarray}
If $E,F$ have the same order, then 
\begin{equation} (\alpha \, E + \beta \, F, 
\gamma \, E + \delta \, F)_{2r+1} = 
(\alpha \, \delta - \beta \, \gamma) \, (E,F)_{2r+1}, 
\label{comb} \end{equation}
for arbitrary constants $\alpha,\beta,\gamma,\delta$. 
This shows that the odd transvectants $(E,F)_{2r+1}$ are 
{\sl combinants} of $E,F$, i.e., up to a scalar, they 
depend only on the pencil spanned by $E,F$. 

If $E,F$ are given {\sl symbolically}, then \cite[\S 49]{GrYo} gives 
an algorithm for calculating their transvectants. See 
Proposition \ref{prop.w} for a typical instance of its use. 

The following lemma is elementary (see \cite[Lemma 2.2]{Goldberg1}). 
\begin{Lemma} \sl Let $E,F$ be nonzero forms of order $e$ such that 
$(E,F)_1=0$. Then $E = \cons F$. \qed 
\label{lemma.jac} \end{Lemma} 

\subsection{} \label{self-dual}
Each representation $S_e$ is self-dual, i.e., 
we have an isomorphism 
\[ S_e \stackrel{\sim}{\lra} S_e^* = \Hom_{SL(V)}(S_e , \complex). \] 
This map sends an order $e$ form $E$ to the functional 
\[\delta_E: S_e \lra \complex, \quad F \lra (E,F)_e. \] 

\subsection{The Gordan series} 
Given three forms $E,F,G$, this very useful series describes certain 
linear dependency relations between transvectants of the 
type $((E,F)_\star,G)_\star$ and $((E,G)_\star,F)_\star$. 

Let $E,F,G$ be of orders $e,f,g$ respectively, and 
$a_1,a_2,a_3$ three integers satisfying the following 
conditions:
\begin{itemize} 
\item 
$\, a_2 + a_3 \le e, \, a_1 + a_3 \le f, \, a_1 + a_2 \le g$, and 
\item either 
$a_1 = 0$ or $a_2 + a_3 = e$ (or both). 
\end{itemize} 
Then we have an identity
\begin{equation} \begin{aligned} 
& \sum\limits_{i \ge 0} \, \frac{\binom{f-a_1-a_3}{i}\binom{a_2}{i}}
{\binom{e+f-2a_3-i+1}{i}} \, ((E,F)_{a_3+i},G)_{a_1+a_2-i} \\ 
= (-1)^{a_1} & \sum\limits_{i \ge 0} \, 
\frac{\binom{g-a_1-a_2}{i}\binom{a_3}{i}}
{\binom{e+g-2a_2-i+1}{i}} \, ((E,G)_{a_2+i},F)_{a_1+a_3-i} \, .  
\label{Gordan} \end{aligned}\end{equation}
Usually (\ref{Gordan}) is denoted by 
$\Gordan{E}{F}{G}{e}{f}{g}{a_1}{a_2}{a_3}$. 

\subsection{The Clebsch-Gordan series} 
Let $\polar$ denote the polarization operator 
\[ y_0 \frac{\partial}{\partial x_0} + y_1 \frac{\partial}{\partial x_1}. \] 
If $E$ is a form of order $e$, then define its $m$-th polar 
to be 
\[ \po{E}{m} = \frac{(e-m)!}{e!} \, (\polar)^m \, E, \] 
which is an $\ux\uy$-form of order $(e-m,m)$. By Euler's theorem, 
we can recover $E$ from $\po{E}{m}$ by the substitution $\uy:=\ux$. 
If $e$ is even, we will denote $\po{E}{e/2}$ by $E^\sympol$. 
It is symmetric in $\ux, \uy$, and naturally thought of as an 
element of $S_2(S_{e/2})$. 

The Clebsch-Gordan series is a more precise version of 
the decomposition (\ref{clebschgordan}). For 
forms $E,F$ of orders $e,f$, it gives an identity 
\begin{equation} E(\ux) \, F(\uy) = \sum\limits_{r \ge 0} \, 
\frac{\binom{e}{r} \binom{f}{r}}{\binom{e+f-r+1}{r}} \, 
\omega^r \, \po{(E,F)_r}{f-r}. 
\label{CG} \end{equation}

\begin{Remark} \rm 
The notional distinction between the Gordan series and 
Clebsch-Gordan series is merely for convenience of reference, and 
has no historical basis. In fact (\ref{CG}) directly leads 
to (\ref{Gordan}) (see \cite[\S 52]{GrYo}). 
\end{Remark} 

Now let $\U \in S_2(S_{d-1})$. We identify $\U$ with 
an $\ux\uy$-form of order $(d-1,d-1)$ which is symmetric in 
both sets of variables. It can then be expressed 
as a `Taylor series' in $\omega$. Define constants 
\begin{equation} 
c_r = \frac{2 \, \binom{d}{2r+1}^2}{\binom{2d-2r}{2r+1}} \quad 
\text{for $0 \le r \le \lfloor \frac{d-1}{2} \rfloor$}. 
\label{cr} \end{equation}
\begin{Proposition} \sl  
There exists a unique sequence of forms 
\[ 
U_\bullet = (U_1,U_3,\dots,U_{2r+1}, \dots), \] 
where $\ord \, U_{2r+1} = 2(d-2r-1)$, such that 
\[ 
\U = \sum\limits_{r \ge 0} c_r \, \omega^{2r} \, (U_{2r+1})^\sympol. \] 
\end{Proposition} 
\demo First we prove the existence.  
Since $\U$ is symmetric in $\ux$ and $\uy$, it is a linear 
combination of expressions of the form 
\[ \langle i \, j \rangle  = 
x_0^{d-1-i} x_1^i \, y_0^{d-1-j}y_1^j + 
x_0^{d-1-j} x_1^j \, y_0^{d-1-i}y_1^i. \] 
Let $A = x_0^{d-1-i}x_1^i, B = x_0^{d-1-j}x_1^j$, so 
\[ \langle i \, j \rangle  = A(\ux)B(\uy) + B(\ux)A(\uy). \] 
Rewrite the right-hand side as a sum of two Clebsch-Gordan series. 
By property (\ref{trans.1}), only the 
even powers of $\omega$ will survive. This shows the existence claim 
for $\langle i j \rangle$, and hence in general by linearity. 

Conversely, let $U_\bullet, U_\bullet'$ be two such sequences for $\U$. By 
the substitution $\uy:=\ux$, we deduce $U_1 = U_1'$. Now 
divide $\U - U_1$ by $\omega^2$ and again let $\uy:=\ux$ etc., 
then we successively see that $U_{2r+1} = U_{2r+1}'$ for all $r$. \qed 

\smallskip 

Henceforth, $A,B$ will always denote linearly independent 
forms of order $d$. We will write 
\begin{equation} T_i:= (A,B)_i, \quad 
\Pi_{A,B}:= \text{Span}\{A,B\} \subseteq S_d. 
\label{T_and_Pi} \end{equation}

\begin{Proposition} \sl 
With notation as above, 
\begin{equation}
\R(A,B) = \sum\limits_{r \ge 0} \, c_r \,\omega^{2r} \, 
(T_{2r+1})^\sympol. \label{series.taylor} \end{equation}
\end{Proposition}
\demo Express $A(\ux)B(\uy)$ and $B(\ux)A(\uy)$ as Gordan 
series and subtract. By property (\ref{trans.1}), only the odd powers 
of $\omega$ will survive. Now divide by $\omega$, then they all 
become even powers. \qed 

\smallskip 

It follows that the collection $\{ T_{2r+1}: r \ge 0\}$ determines 
$\R(A,B)$. It will be shown below that the terms $r=0,1$ are already 
sufficient. 
\section{The Wronskian o.d.e.} \label{section.wr.ode}
\subsection{Generalities on Wronskians} 
Given integers $p,q$ with $q \le p+1$, there is an 
isomorphism of $SL(V)$-modules (see \cite[\S 11]{FH})
\begin{equation} \wedge^q S_p \stackrel{\sim}{\lra} S_q(S_{p-q+1}). 
\label{extsym} \end{equation}
Composing it with the natural surjection 
\begin{equation} S_q(S_{p-q+1}) \lra S_{q(p-q+1)}, 
\label{surj}\end{equation}
we get the Wronskian map 
\[ \Theta: \wedge^q S_p \lra S_{q(p-q+1)}. \] 
If $F_1, \dots, F_q$ are order $p$ forms, then 
their Wronskian 
$\Theta(F_1 \wedge \dots \wedge F_q)$ is given 
by the $q \times q$ determinant 
\[ (i,j) \lra 
\frac{\partial^{\,q-1} F_i}{\partial x_0^{q-j} \, \partial x_1^{j-1}}
\quad 1 \le i,j \le q. \] 
The crucial property of the construction is that 
$\Theta$ is nonzero on decomposable tensors, i.e., 
$\Theta(F_1 \wedge \dots \wedge F_q)=0 \iff 
F_1 \wedge \dots \wedge F_q=0 \iff $
the $F_i$ are linearly dependent.

\subsection{} 
Now let $A,B,F$ be of order $d$, with Wronskian 
\[ \W = \Theta(A \wedge B \wedge F) = 
 \left| \begin{array}{ccc} 
A_{x_0^2} & A_{x_0x_1} & A_{x_1^2} \\ 
B_{x_0^2} & B_{x_0x_1} & B_{x_1^2} \\ 
F_{x_0^2} & F_{x_0x_1} & F_{x_1^2}
\end{array} \right|. \] 
We will evaluate $\W$ symbolically. Let us write 
\begin{equation} A = \alpha_x^d, \; B = \beta_x^d, \; F = f_x^d. 
\label{sym.exp} \end{equation}
As usual, $\alpha_x$ stands for the symbolic linear form 
$\alpha_0 x_0 + \alpha_1 x_1$, and 
$\br{\alpha}{\beta}$ for $\alpha_0 \beta_1 - \alpha_1 \beta_0$ etc. 

\begin{Lemma} \sl With notation as above, 
\begin{equation} 
\frac{1}{(d^2-d)^3} \, 
\W = \br{\alpha}{\beta}\br{\alpha}{f}\br{\beta}{f} \, 
\alpha_x^{d-2} \, \beta_x^{\, d-2} \, f_x^{\, d-2}. 
\label{sym.wronskian} \end{equation} 
\end{Lemma} 
\demo Differentiating (\ref{sym.exp}), we get expressions such as 
\[ A_{x_0x_1} = d(d-1) \, \alpha_x^{d-2} \, \alpha_0 \, \alpha_1. \]
Substitute these into $\W$ and factor out  
$\alpha_x^{d-2} \, \beta_x^{d-2} \, f_x^{d-2}$. We are left 
with a Vandermonde determinant which evaluates to 
$\br{\alpha}{\beta}\br{\alpha}{f}\br{\beta}{f}$. \qed 

Now we will rewrite $\W$ in terms of transvectants. 
\begin{Proposition} \sl 
With notation as in (\ref{T_and_Pi}), we have an identity 
\[ \frac{1}{(d^2-d)^3} \, \W = (T_1,F)_2 - \frac{d-2}{4d-6} \, T_3 \, F. \] 
\label{prop.w} \end{Proposition} 
\demo 
Symbolically, the transvectants can be written as 
\[ 
T_1 = \br{\alpha}{\beta} \, \alpha_x^{\,d-1} \beta_x^{\,d-1}, \quad 
T_3 = \br{\alpha}{\beta}^3 \, \alpha_x^{\,d-3} \beta_x^{\,d-3}. \] 
First we calculate the transvectant $(T_1,F)_2$ using the 
algorithm given in (see \cite[\S 49]{GrYo}): 
\begin{itemize} 
\item 
Calculate the second polar $T_1$. It is equal to 
\[ \begin{aligned} 
{} & \frac{(2d-4)!}{(2d-2)!} \, (\polar)^2 \, T_1  = 
\frac{1}{(2d-2)(2d-3)} \br{\alpha}{\beta} \, \alpha_x^{\,d-3} \beta_x^{\,d-3} 
\times \\ 
& \{(d-1)(d-2) \, \alpha_x^2 \, \beta_y^2 + 
 2 (d-1)^2 \, \alpha_x \, \beta_x \, \alpha_y \, \beta_y + 
 (d-1)(d-2) \, \beta_x^2 \, \alpha_y^2 \}. 
\end{aligned} \] 
\item
Make substitutions $\alpha_y:=\br{\alpha}{f},\beta_y:=\br{\beta}{f}$, 
and multiply by $f_x^{\,d-2}$. The result is 
\begin{equation} \begin{aligned} 
(T_1,F)_2 = &  
\frac{d-2}{4d-6} \br{\alpha}{\beta} \, 
\alpha_x^{\,d-3} \beta_x^{\,d-3} f_x^{\,d-2} \times  \\
& \{ \br{\beta}{f}^2 \, \alpha_x^2 + \frac{2d-2}{d-2} 
\br{\alpha}{f} \br{\beta}{f} \, 
\alpha_x \beta_x + \br{\alpha}{f}^2 \, \beta_x^2 \}. 
\end{aligned} \label{t1f2} \end{equation}
\end{itemize} 

We would like to compare (\ref{sym.wronskian}) and (\ref{t1f2}), 
so we will rewrite both of them in terms of 
{\sl standard monomials} (see \cite[Ch.~3]{Sturmfels}). Order the variables as 
$\alpha < \beta < f < x$. 
The monomial $\br{\beta}{f} \, \alpha_x$ is nonstandard, so 
use the Pl{\"u}cker syzygy to rewrite it as 
\[ \br{\beta}{f} \, \alpha_x = 
\br{\alpha}{f} \, \beta_x - \br{\alpha}{\beta} f_x. \] 
Substitute this into the right hand sides of (\ref{sym.wronskian}) and 
(\ref{t1f2}). Subtracting the two expressions, we get 
\[ (T_1,F)_2 - \frac{1}{(d^2-d)^3} \W = 
\frac{d-2}{4d-6} \, \br{\alpha}{\beta}^3 \, \alpha_x^{d-3}\beta_x^{d-3}f_x^d 
= \frac{d-2}{4d-6} \, T_3 \, F. \] 
This completes the proof. \qed 

\subsection{} 
If $M,N$ are forms of orders $2d-2,2d-6$ respectively, then 
we define 
\begin{equation} 
\psi_{M,N}(F) := (M,F)_2 - \frac{d-2}{4d-6} \, N \, F. 
\label{psi} \end{equation}
For fixed $M,N$, we are interested in the differential equation 
\begin{equation} \psi_{M,N}(F) = 0. \label{diffeq} \end{equation} 
We may call this the Wronskian (second order) ordinary differential 
equation with parameters $M,N$. (It is always assumed that $M \neq 0$, 
otherwise the equation is of no interest.) 
The following corollary is immediate. 
\begin{Corollary} \sl If $F$ is of order $d$, then 
$F \in \Pi_{A,B}$ iff $\psi_{T_1,T_3}(F) = 0$. 
\label{cor.T1T3} \end{Corollary} 
\demo Indeed, $\psi_{T_1,T_3}(F)=0$ iff $A,B,F$ are linearly 
dependent. \qed 

\smallskip 

Hence, given $T_1,T_3$, the pair $\{A,B\}$ is determined up to a 
unimodular transformation (cf.~(\ref{comb})). It follows that $T_1,T_3$ 
together determine all the $T_{2r+1}$. 

\begin{Proposition} \sl 
Let $M,N$ be of orders $2d-2,2d-6$. Assume that (\ref{diffeq}) has 
two linearly independent solutions $A,B$ of order $d$. Then there 
exists a nonzero constant $\lambda$ such that 
$M = \lambda \, T_1, \; N = \lambda \, T_3$. 
\label{prop.pencil} \end{Proposition} 
\demo 
Multiply the identities $\psi_{M,N}(A)=0,\psi_{M,N}(B)=0$ 
by $B,A$ respectively and subtract, this gives 
$ B\, (M,A)_2 =  A \, (M,B)_2$. Now the Gordan series 
\[ \Gordan{A}{M}{B}{d}{2d-2}{d}{0}{0}{2}, \quad 
   \Gordan{B}{M}{A}{d}{2d-2}{d}{0}{0}{2} \] 
respectively give identities 
\[ \begin{aligned} 
(A,M)_2 \, B & = (AB,M)_2 + ((A,B)_1,M)_1 + \frac{d}{4d-2} (A,B)_2 \,M \\ 
(B,M)_2 \, A & = (BA,M)_2 + ((B,A)_1,M)_1 + \frac{d}{4d-2} (B,A)_2 \,M.
\end{aligned} \] 
Subtracting and using property (\ref{trans.1}) for $A,B$, we get 
$((A,B)_1,M)_1 =0$. Now $(A,B)_1 \neq 0$ since $A,B$ are 
independent, but then Lemma \ref{lemma.jac} implies 
that $M = \lambda \, (A,B)_1$ for some $\lambda$. Finally 
\[ \frac{d-2}{4d-6} \, N A = (M,A)_2 = \lambda \, (T_1,A)_2 = 
\lambda \, \frac{d-2}{4d-6} \, T_3 \, A, \] 
hence $N = \lambda \, T_3$. \qed 

\subsection{}
We have shown that the following conditions are equivalent for 
the pair $(M,N)$. 
\begin{itemize} 
\item[(i)] 
There exist $A,B$ such that 
$M = (A,B)_1, N = (A,B)_3$. 
\item[(ii)]
There exist $A,B$ such that 
\[ \R(A,B) = c_0 \, M^\sympol + c_1 \, \omega^2 \, N^\sympol + O(\omega^4). \] 
\item[(iii)] 
The dimension of the kernel of the map 
\[ \psi_{M,N}: S_d \lra S_{3d-6} \] 
is at least $2$ (and then it is exactly $2$).
\end{itemize} 
We can now construct the covariants $\Phi_r$ as in the 
introduction. Clearly (iii) is equivalent to the condition that the map 
\[ \wedge^d \psi_{M,N}: 
\wedge^d S_d \lra \wedge^d (S_{3d-6}) \] 
be zero. Identify $\wedge^d S_d$ with $S_d$ via (\ref{extsym}). 
Let $f_1$ denote the image of $\wedge^d \psi_{M,N}$ via the 
isomorphism 
\[ \Hom_{SL(V)}(S_d, \wedge^d S_{3d-6}) \simeq 
\Hom_{SL(V)}(\complex, \wedge^d S_{3d-6} \otimes S_d).\]
Consider the composite morphism 
\[ \complex \stackrel{f_1}{\lra} 
\wedge^d S_{3d-6} \otimes S_d \stackrel{f_2}{\lra} 
S_d (S_{2d-5}) \otimes S_d \stackrel{f_3}{\lra} 
S_{d(2d-5)} \otimes S_d, \] 
where $f_2$ comes from the isomorphism (\ref{extsym}), and 
$f_3$ from the natural surjection (\ref{surj}). For each 
$0 \le r \le d$, we have projection maps 
\[ \pi_r: S_{d(2d-5)} \otimes S_d \lra S_{d(2d-4)-2r} \] 
induced by the decomposition (\ref{clebschgordan}). 

Define $\Phi_r(M,N)$ to be the image of $1 \in \complex$ 
via the map $\pi_r \circ f_3 \circ f_2 \circ f_1$. This is 
a joint covariant of $M,N$ of order $d(2d-4)-2r$. 
We will describe it in coordinates. 
For $0 \le i \le d$, define 
\[ w_i = (-1)^i \frac{1}{i!(d-i)!} \, 
\Theta( \, \bigwedge\limits_{\stackrel{s=0}{s \neq i}}^d 
\psi_{M,N}(x_0^s \, x_1^{d-s})), \] 
which is an element of $S_{d(2d-5)}$. Then 
\[ (f_3 \circ f_2 \circ f_1) (1) = 
\sum\limits_{i=0}^d w_i \otimes x_0^i \, x_1^{d-i}, \quad 
\text{and} \quad 
\Phi_r = \sum\limits_{i=0}^d \, (w_i, x_0^i \, x_1^{d-i})_r. \] 
All of this is straightforward and follows by chasing through the 
$f_i$. Each $\Phi_r$ has total degree $d$ in the coefficients 
of $M,N$ (because $w_i$ does). 

\begin{Theorem} \sl 
Let $M,N$ be orders $2d-2,2d-6$ respectively. 
Then the pair $(M,N)$ satisfies the (equivalent) conditions 
(i)--(iii) if and only if 
\[ \Phi_r(M,N) =0 \quad \text{for $0 \le r \le d$.} \] 
\end{Theorem} 
\demo If (iii) holds, then $f_1=0$, which  shows the 
`only if' part. Conversely, assume that all 
the $\Phi_r$ vanish. Then $(f_3 \circ f_2 \circ f_1)(1) =0$, which 
implies that all the $w_i$ vanish. By the fundamental property 
of Wronskians, the forms 
\[ \psi_{M,N}(x_0^s \, x_1^{d-s}), \quad 0 \le s \le d, \, s \neq i \] 
are linearly dependent for any $i$. But then the map 
$\wedge^d \psi_{M,N}$ is zero on every basis element of $\wedge^d S_d$, 
hence it is zero. This implies (iii). 
\qed 

\subsection{The incomplete Pl{\"u}cker imbedding} \label{incomplete.plucker} 
The fact that $\R$ is determined by $T_1,T_3$ has the following geometric 
interpretation. Assume $d \ge 3$, and let 
$\G = G(2,S_d)$ denote the Grassmannian of two-dimensional 
subspaces in $S_d$. (See \cite[Lecture 6]{JoeH} for generalities on 
Grassmannians.) The line bundle $\O_{\G}(1)$ has 
global sections 
\[ H^0(\G,\O_\G(1)) \simeq \wedge^2 S_d \simeq S_2(S_{d-1}) 
\simeq \bigoplus\limits_{r \ge 0} S_{2d-4r-2}. \] 
The usual Pl{\"u}cker imbedding is given by the complete 
linear system $|\O_\G(1)|$. 
Consider the subspace $W = S_{2d-2} \oplus S_{2d-6} 
\subseteq H^0(\O_\G(1))$. 

\begin{Proposition} \sl 
The map 
\[ \mu : \G \lra \P \, W, \quad 
\P \, \Pi_{A,B} \lra [ \, T_1 \oplus T_3 \, ] \] 
is an isomorphic imbedding. 
\end{Proposition} 
The usual conventions (\cite[Ch.~II, \S 7]{Ha}) dictate that the 
imbedding is in $\P W^*$, but note the self-duality in \S \ref{self-dual}.

\smallskip 

\demo We have already shown that $\mu$ is a set-theoretic injection. 
To complete the proof, it suffices to show that it is an injection 
on tangent spaces at every point (cf.~\cite[Ch.~II, Prop.~7.3]{Ha}). 
The Zariski tangent space to $\G$ at $\Pi = \Pi_{A,B}$ is 
canonically isomorphic to $\Hom(\Pi,S_d/\Pi)$ 
(see \cite[Lecture 16]{JoeH}). 
Let $\alpha: \Pi \lra S_d/\Pi$ be a tangent vector, and say 
\[ \alpha(A) = Q + \Pi, \quad \alpha(B) = P + \Pi, \] 
for some forms $P,Q$ of order $d$. 

The tangent space to $\P W $ at $[T_1 \oplus T_3]$ is 
isomorphic to $W/\langle T_1 \oplus T_3 \rangle$. 
Let $d \mu: T_{\G,\Pi} \lra T_{\P,\mu(\Pi)}$ denote the induced map 
on tangent spaces. Then $d\mu(\alpha)$ is the element 
\[ ((A,P)_1+(Q,B)_1) \oplus ((A,P)_3+(Q,B)_3) \in W \] 
considered modulo $T_1 \oplus T_3$. (To see this, 
let $\epsilon$ be an `infinitesimal'. Now 
expand $(A + \epsilon \, Q, B + \epsilon \, P)_i, i=1,3$, and 
set $\epsilon^2 = 0$.) 
 
We would like to show that $d\mu$ is injective, hence 
suppose that $d\mu(\alpha)=0$. Then there exits a constant $c$ such that 
\[ \begin{aligned} 
(A,P)_1 + (Q,B)_1 & = c \, (A,B)_1, \\
(A,P)_3 + (Q,B)_3 & = c \, (A,B)_3. \end{aligned} \] 
Substitute $P + c \, B$ for $P$ (which does not change $\alpha$), then 
\[ (A,P)_1 = (B,Q)_1, \quad (A,P)_3 = (B,Q)_3. \] 
If the first pair is zero, then $P,Q$ are respectively 
constant multiples of $A,B$, hence $\alpha=0$. If not, then 
$\Pi_{A,P} = \Pi_{Q,B}$ by Corollary \ref{cor.T1T3}. But this implies 
$\Pi_{A,B} = \Pi_{P,Q}$, again forcing $\alpha=0$. \qed 

\begin{Remark} \rm 
Let ${\mathfrak a} \subseteq \text{Sym}^\bullet \, W$ denote the ideal 
generated by the coefficients of $\Phi_0,\dots,\Phi_d$, and 
$J$ the homogeneous ideal of the image $\mu(\G) \subseteq \P W$. 
Since ${\mathfrak a}$ defines the image set-theoretically, 
$(\sqrt{\mathfrak a})_{\text{sat}} = J$. Already for $d=3$
these ideals differ (since ${\mathfrak a}$ is generated in 
degree $3$ and $J$ in degree $2$), but I do not know if one can state 
a more precise relation in general. 
\end{Remark} 

\section{Generic reduction formulae} \label{section.gen.red}
\subsection{} 
We begin with the example which eventually led to 
the main result of this section. Let $A,B$ be of order $2$. 
The series $\Gordan{A}{B}{A}{2}{2}{2}{0}{1}{1}$ implies the relation
\[ ((A,B)_1,A)_1 + \frac{1}{2} \, (A,B)_2 \, A = 
\frac{1}{2} \, (A,A)_2 \, B; \] 
which can be rewritten as 
\[ -\frac{2}{(A,A)_2} (A,T_1)_1 = B - \frac{(A,B)_2}{(A,A)_2} \, A. 
\] 
Hence, given $\R$ (which involves only $T_1$ in this case) and $A$, 
the function 
\begin{equation} (A,T_1) \lra  -\frac{2}{(A,A)_2} (A,T_1)_1 
\label{red1} \end{equation}
recovers $B$ up to an additive multiple of $A$. 
(Since $\R(A,B + \cons A) = \R(A,B)$, the last proviso is 
indispensable.) 
We will show that there exist such formulae for every $d$. 
\begin{Remark} \rm 
We may call (\ref{red1}) a reduction formula in the following 
sense. If we are given a linear second order 
o.d.e., together with one of its solutions, then a second solution 
can be found by the method of `reduction of order' 
(see \cite[\S 44]{Rainville}). In our case, we are to find $B$, 
given the equation $\psi_{T_1,T_3}(F)=0$ with one solution $A$. 
However, this analogy is inexact in two respects: 
\begin{itemize} 
\item 
our formula will involve all the $\{T_{2r+1}\}$, and 
not merely $T_1,T_3$, 
\item 
the process is algebraic and involves no integration. 
\end{itemize} 
Moreover, the formula is generic in the sense that it is only defined over 
an open subset, e.g., the set $\{A \in \P^2: (A,A)_2 \neq 0\}$ above. 
\end{Remark}

\subsection{} 
Throughout this section we assume that $A,B$ are order $d$ forms 
whose coefficients are algebraically independent indeterminates. 
Write 
\begin{equation} 
A = \sum\limits_{p=0}^d \binom{d}{p} \, a_p \, x_0^{d-p} \, x_1^p. 
\label{exp.A} \end{equation}
Let $\J$ be an invariant of $A$ of degree (say) $n$. We define its 
first evectant (cf.~\cite{Sylvester3}) to be 
\begin{equation}
\Ev_\J = \frac{(-1)^d}{n} \sum\limits_{q=0}^d \, 
(-1)^q \, \frac{\partial \J}{\partial a_q} \, x_0^q \, x_1^{d-q},
\label{Ev} \end{equation}
it is a covariant of degree-order $(n-1,d)$. The initial scaling 
factor is chosen so as to make the following lemma true: 
\begin{Lemma} \sl 
We have an identity $(\Ev_\J,A)_d  = \J$. 
\label{EV} \end{Lemma} 
\demo 
Substitute (\ref{exp.A}) and (\ref{Ev}) in formula (\ref{trans.defn}). 
We get a nonzero term whenever $p=q$ and $i=d-p$, hence 
\[ \begin{aligned} (\Ev_\J,A)_d  & = 
\frac{(-1)^{2d}}{n \, (d \, !)^2} \, \sum\limits_{p=0}^d \, 
(p!(d-p)!)^2 \binom{d}{d-p}^2 \, a_p \, \frac{\partial \J}{\partial a_p} \\ 
& = \frac{1}{n} \sum\limits_p  \, a_p \, \frac{\partial \J}{\partial a_p}
= \J,  \end{aligned} \] 
the last equality is by Euler's theorem. \qed 

\smallskip 

Now our generic reduction formula is as follows. Let 
\begin{equation} 
\beta(A,\R) = - \frac{1}{\J} 
\sum\limits_{r \ge 0} \, c_r \, (\Ev_\J,T_{2r+1})_{d-2r-1},
\end{equation}
with the $c_r$ as in (\ref{cr}). 
\begin{Theorem} \sl 
With notation as above, 
\[ \beta(A,\R) = B - \frac{(\Ev_\J,B)_d}{\J} \, A. \] 
\label{Theorem.red.formula} \end{Theorem}
Hence, as long as $A$ stays away from the hypersurface $\{ \J =0\}$, 
we can recover $B$ from $A$ and $\R(A,B)$. 

\begin{Remark} \rm 
If $d$ is even, then we can take $\J$ to be the unique 
degree two invariant $(A,A)_2$. 
There is no invariant in degrees $\le 3$ if $d$ is odd, but then 
there exists a degree four invariant $\J = ((A,A)_{d-1},(A,A)_{d-1})_2$. 
\end{Remark} 

\subsection{} 
The proof of the theorem will emerge from the discussion below. 
The element $A\wedge B \in \wedge^2 S_d$ defines a map 
\[ \sigma_{A \wedge B}: S_d \lra S_d, \quad 
F \lra (F,B)_d \, A - (F,A)_d \, B. \] 
We identify the codomain of $\sigma=\sigma_{A \wedge B}$ with 
$S_d^*$ as in \S \ref{self-dual}. 

\begin{Lemma} \sl With the convention above, $\sigma$ 
is skew-symmetric, i.e., 
\[ 
\delta_{\sigma(F)}(G) = - \delta_{\sigma(G)}(F), 
\quad \text{for $F,G \in S_d$.} \] 
\end{Lemma} 
\demo 
Unwinding the definitions, this becomes 
\[ (F,B)_d (A,G)_d - (F,A)_d (B,G)_d = 
- \{(G,B)_d (A,F)_d - (G,A)_d(B,F)_d \}, \] 
which is obvious. \qed 

\begin{Lemma} \sl 
With notation as above, 
\begin{equation} 
\sigma(F) = [ \, (F,\R)_{d-1}^\ux \, ]_{\uy:=\ux}. 
\label{sigma.F.R} \end{equation} 
\end{Lemma}
The right hand side of this equation is interpreted as follows: 
calculate the $(d-1)$-th transvectant of $F$ and $\R$ as 
$\ux$-forms (treating the $\uy$ in $\R$ as constants). This produces 
an $\ux\uy$-form of order $(1,d-1)$; finally replacing $\uy$ by 
$\ux$ gives a form of order $d$. 

\smallskip 

\demo 
We will calculate both sides symbolically. Let 
$A = \alpha_x^d, \; B = \beta_x^d, \; F = f_x^d$, then 
\[ \begin{aligned} 
\R & = \frac{A(\ux)B(\uy) - A(\uy)B(\ux)}{\omega} 
= \frac{\alpha_x^d \, \beta_y^d - \alpha_y^d \, \beta_x^d}{\omega} \\ 
& = \frac{
(\alpha_x \, \beta_y - \alpha_y \, \beta_x)
\sum\limits_{i=0}^{d-1} \, 
(\alpha_x \, \beta_y)^{d-1-i}(\alpha_y \, \beta_x)^i}
{\omega} \\ 
& = \br{\alpha}{\beta} \, 
\sum\limits_i (\alpha_x \, \beta_y)^{d-1-i}(\alpha_y \, \beta_x)^i.
\end{aligned} \] 

Now calculate the $(d-1)$-th transvectant of $F$ 
with each summand in the last expression. (We have agreed to 
treat $\alpha_y,\beta_y$ as constants). Using the 
algorithm of \cite[\S 49]{GrYo}, 
\[ (f_x^{\,d},\alpha_x^{d-1-i} \, \beta_x^i)_{d-1} = 
(-1)^{d-1} \br{\alpha}{f}^{d-1-i} \, \br{\beta}{f}^i \, f_x. \] 
Hence, 
\begin{equation} 
[\, (F,\R)_{d-1}^\ux \, ]_{\uy:=\ux} = 
(-1)^{d-1} \br{\alpha}{\beta} \, f_x \sum\limits_i \, 
\br{\alpha}{f}^{d-1-i} \, \br{\beta}{f}^i \, 
\alpha_x^i \, \beta_x^{d-1-i}. 
\label{FR} \end{equation} 
Now directly from the definition, 
\[ \begin{aligned} \sigma(F) & =  
\{ \br{f}{\beta}^d \, \alpha_x^d - \br{f}{\alpha}^d \, \beta_x^d \} \\
& = (-1)^d 
\{ \br{\beta}{f}^d \, \alpha_x^d - \br{\alpha}{f}^d \, \beta_x^d \} \\
& = (-1)^d \{\br{\beta}{f} \, \alpha_x - \br{\alpha}{f} \, \beta_x \} 
\sum\limits_i \, \br{\alpha}{f}^{d-1-i} \, \br{\beta}{f}^i \, 
\alpha_x^i \, \beta_x^{d-1-i}. 
\end{aligned} \] 
Since $\br{\beta}{f}\alpha_x - \br{\alpha}{f}\beta_x 
= -\br{\alpha}{\beta}f_x$, the last 
expression is identical to (\ref{FR}). \qed 

\smallskip 

\begin{Lemma} \sl 
Let ${\mathbf T}$ be an arbitrary form of order $2d-4r-2$. Then  
\[ [\, (F,\omega^{2r} \, T^\sympol)_{d-1}^\ux \,]_{\uy:=\ux} 
= (F,{\mathbf T})_{d-2r-1}. \] 
\end{Lemma} 
\demo Let ${\mathbf T} = t_x^{\, 2d-4r-2}$, so that 
$T^\sympol = t_x^{\, d-2r-1} \, t_y^{\, d-2r-1}$. Then make a 
calculation as in the previous lemma. \qed 

\smallskip 

Now substitute the Taylor series (\ref{series.taylor}) into the 
right hand side of (\ref{sigma.F.R}), and use the previous lemma. 
This gives the formula 
\begin{equation} 
\sigma(F) = \sum\limits_{r \ge 0} \, c_r \, (F,T_{2r+1})_{d-2r-1}. 
\label{sigma_and_T} \end{equation}
Now specialize to $F = \Ev_\J$. Then 
\[ 
\sigma(\Ev_\J) = (\Ev_\J,B)_d \, A - (\Ev_\J,A)_d \, B = 
(\Ev_\J,B)_d \, A - \J \, B, \] 
hence 
\[ \beta(A,\R) = 
-\frac{1}{\J} \, \sigma(\Ev_\J) = B - \frac{(\Ev_\J,B)_d}{\J} \, A. \] 
This completes the proof of Theorem \ref{Theorem.red.formula}. \qed 

\section{Formulae for $T_5$ and $T_7$} \label{section.t5t7}
We have observed that $T_1,T_3$ determine 
the higher odd transvectants $T_{2r+1}$. However this dependence is 
rather indirect, and it is unclear if one can give a formula for the latter 
in terms of the former. In this section we give such 
explicit formulae for $T_5$ and $T_7$. 
\subsection{The Pl{\"u}cker relations} 
Let 
\[ \G \subseteq \P (\wedge^2 S_d) = 
\P ( \, \bigoplus\limits_{r \ge 0} \, S_{2d-4r-2}) 
\] 
be the usual Pl{\"u}cker imbedding, and let $\I$ denote the 
homogeneous ideal of the image. It is well-known that
$\I$ is generated by its quadratic part $\I_2$, 
usually called the module of Pl{\"u}cker relations. 
\begin{Lemma} \sl 
As $SL(V)$-modules, $\I_2 \simeq \wedge^4 S_d$. 
\end{Lemma} 
\demo 
Consider the short exact sequence 
\[ 
0 \ra \I_2 \ra H^0(\O_{\P(\wedge^2 S_d)}(2)) \ra 
H^0(\O_{\G}(2)) \ra 0. \] 
(The exactness on the right comes from the projective normality of 
the imbedding.) Using the plethysm formula of 
\cite[\S I.8, Example 9]{MacDonald}, the middle term is 
isomorphic to 
\[ S_2(\wedge^2 S_d) \simeq S_{(2,2)} (S_d) \oplus \wedge^4  S_d. \] 
By the Borel-Weil theorem (see \cite[p.~687]{Porras}), 
$H^0(\O_\G(2)) \simeq S_{(2,2)}(S_d)$. This completes the 
proof. \qed 

\smallskip 

Each Pl{\"u}cker relation corresponds to an algebraic identity 
between the $\{T_{2r+1}\}$. To make this more 
precise, let $\{M_{2r+1}:r \ge 0\}$ be generic forms of orders $2d-4r-2$, and 
$S_e \stackrel{\xi}{\hookrightarrow} \I_2$ an inclusion of 
$SL(V)$-modules. Then $\xi$ corresponds to a joint covariant 
$\Xi \, (M_1,M_3,\dots)$ of order $e$ and total degree two in the 
$\{M_{2r+1} \}$, such that 
\[ \Xi \, (T_1,T_3,\dots) = 0, \qquad \text{for any $A,B$ of order $d$.} 
\] 

\begin{Example} \rm 
Assume $d=4$. In this case $\I_2 \simeq S_4$, so we look for 
an order $4$ covariant in $M_1,M_3$. There are 
three `monomials' of total degree $2$ and order $4$, namely 
$(M_1,M_1)_4,(M_1,M_3)_2,M_3^2$. Our covariant must be a linear 
combination of these, i.e., 
\[ \Xi(M_1,M_3) = \alpha_1 \, (M_1,M_1)_4 + 
\alpha_2 \, (M_1,M_3)_2 + \alpha_3 \, M_3^2, \] 
for some constants $\alpha_i$. 

Now specialize to $A = x_0^4,B = x_1^4$, and use formula 
(\ref{trans.defn}) to calculate $T_1,T_3$ and $\Xi$ explicitly. 
Since $\Xi(T_1,T_3)$ must vanish identically, its coefficients 
give $5$ linear equations for the $\alpha_i$. 
Solving these (they must admit a nontrivial solution), we deduce that 
\[ [\alpha_1:\alpha_2:\alpha_3] = [25:-10:-4], \] 
which determines $\Xi$ (of course, up to a scalar). This 
`method of undetermined coefficients' (specializing the 
forms followed by solving linear equations) will be liberally 
used in the sequel. 
\end{Example} 
\begin{Example} \rm 
For $d=3$, the Grassmannian is a quadric hypersurface 
defined by 
\[ \Xi(M_1,M_3) = (M_1,M_1)_4 - \frac{1}{6} \, M_3^2. \] 
\end{Example} 

\subsection{} We begin with a technical 
lemma about the irreducible submodules of $\I_2$. 
\begin{Lemma} \sl If $d \ge 4$, then there exists exactly one copy 
each of the modules $S_{4d-12},S_{4d-16}$ inside $\I_2$.  
\end{Lemma} 
\demo There are isomorphisms 
\[ \I_2 \simeq \wedge^4 S_d \simeq S_4(S_{d-3}) \simeq S_{d-3}(S_4), \] 
where the second isomorphism is from (\ref{extsym}), and the 
third is Hermite reciprocity. Hence we may as well work with 
$S_{d-3}(S_4)$. Now the following are in bijective correspondence 
(see \cite{Littlewood1} for details): 
\begin{itemize} 
\item inclusions $S_e \subseteq S_{d-3}(S_4)$ of $SL(V)$-modules, 
\item 
covariants of degree-order $(d-3,e)$ (distinguished up to scalars) 
for binary quartics, 
\end{itemize} 
Fortunately, a complete set of generators for the covariants of 
binary quartics is known (see \cite[\S 89]{GrYo}). 
It contains five elements, conventionally called $f,H,t,i,j$, having 
degree-orders 
\[ (1,4), \; (2,4), \; (3,6), \; (2,0), \; (4,0). \] 
(It is unnecessary for us to know how they are defined.) 
Each covariant of quartics is a polynomial in the elements of this set. 

Now it is elementary to see that only one expression of 
degree-order $(d-3,4d-12)$ is possible, namely $f^{d-3}$. 
Similarly, the only possible expression for degree-order 
$(d-3,4d-16)$ is $f^{d-5} \, H$. Hence there is exactly 
one copy each of $S_{4d-12}$ and $S_{4d-16}$. \qed 

\subsection{} 
We will find the joint covariant $\Xi$ corresponding to 
$S_{4d-12} \subseteq \I_2$. We look for degree two monomials 
of order $4d-12$ in the $\{T_{2r+1}\}$; any such monomial must be 
of the form 
\[(T_{2a+1},T_{2b+1})_s, \] 
where 
\begin{itemize} 
\item 
$(2d-4a-2)+(2d-4b-2) -2s = 4d-12$, 
\item 
$a,b \le \lfloor \frac{d-1}{2} \rfloor$, 
\item 
$s \le \min \{2d-4a-2,2d-4b-2 \}$, and 
\item 
if $a=b$, then $s$ is even. 
\end{itemize} 
The first condition comes from the order, the rest are forced by 
properties (\ref{trans.1}),(\ref{trans.2}) of transvectants. 
Sifting through these conditions gives only four possibilities, namely
\[ (T_1,T_1)_4, \; (T_1,T_3)_2, \; 
   T_3^2, \; T_1T_5. \] 
Hence we have an identical relation of the form 
\[ \alpha_1 \, (T_1,T_1)_4 + \alpha_2 (T_1,T_3)_2 + 
   \alpha_3 \, T_3^2 - \alpha_4 \, T_1 \, T_5 =0. \] 
Specialize $A,B$ successively to the pairs 
\[ (x_0^d,x_1^d), \; (x_0^{d-1}x_1,x_1^d), \; 
(x_0^{d-2}x_1^2,x_1^d), \; (x_0^{d-1}x_1,x_0 \, x_1^{d-1}), \; \] 
and use the method of undetermined coefficients. 
Up to a scalar, the solution is 
\[ \begin{array}{ll}
\alpha_1 = -\frac{2(2d-3)^2}{d(d-2)} & 
\alpha_2 = \frac{4(2d-3)(d-3)}{d(d-2)} \\ 
\alpha_3 = 1 & 
\alpha_4 = \frac{(d-3)(d-4)(2d-3)^2}{d(2d-5)(2d-7)(d-2)}. 
\end{array} \] 
This gives a formula for $T_5$. 
\begin{Theorem} \sl Assume $d \ge 5$, then 
\[ 
T_5 = \frac{1}{T_1} \, 
( \, \frac{\alpha_1}{\alpha_4} \, (T_1,T_1)_4 
 +\frac{\alpha_2}{\alpha_4} \, (T_1,T_3)_2 
 +\frac{\alpha_3}{\alpha_4} \, T_3^2).  
\] 
\end{Theorem} 
We can make a similar argument with $S_{4d-16}$, this leads to 
a formula for $T_7$. Define 
\[ \begin{array}{ll}
\beta_1 = - \frac{8(2d-5)(2d-7)(2d-3)}{(d-1)(4d-13)} & 
\beta_2 = - \frac{60(2d-7)(2d-5)}{d(d-1)(4d-13)} \\ 
\beta_3 = \frac{12(2d-3)(d-5)}{d(4d-13)} & 
\beta_4 = \frac{20(2d-5)(2d-7)(d-3)}{(d-1)(4d-13)(2d-3)} \\
\beta_5 = 1  & 
\beta_6 = \frac{(d-5)(d-6)(2d-3)(2d-5)}{d(d-1)(2d-9)(2d-11)} 
\end{array} \] 

\begin{Theorem} \sl Assume $d \ge 7$, then 
\[ 
\begin{aligned} 
T_7 = \frac{1}{T_1} 
(& \, \frac{\beta_1}{\beta_6} \, (T_1,T_1)_6 + 
 \frac{\beta_2}{\beta_6} \, (T_1,T_3)_4 +  
 \frac{\beta_3}{\beta_6} \, (T_1,T_5)_2 +  \\
 & \frac{\beta_4}{\beta_6} \, (T_3,T_3)_2 +  
 \frac{\beta_5}{\beta_6} \, T_3 \, T_5). 
\end{aligned} \] 
\end{Theorem} 
(Of course we can substitute for $T_5$ using the previous 
result, but this would make the formula very untidy.) 

This method breaks down for higher transvectants, so a new idea 
will be needed for the general case. 
My colleage A.~Abdesselam, when shown the formulae above, 
remarked that the coefficients look very similar to those appearing in the 
classical hypergeometric series. Perhaps there is something 
to this suggestion. 

\section{Open problems} \label{section.open}
This section contains a series of miscellaneous calculations 
and examples, all of them for small specific values of $d$. 
They should serve simultaneously as a source of 
open questions and further lines of inquiry. 
\subsection{The Jacobian predicate}
Let $A,M$ be forms of orders $d,2d-2$. Consider the 
following predicate 
\[ J(A,M): \text{there exists an order $d$ form $B$ such that 
$(A,B)_1 = M$.} \] 
If $J(A,M)$ holds, then 
$(A,M)_2 = \cons \, T_3 \, A$, hence $A$ must divide $(A,M)_2$. 
We will see below that this condition is sufficient for $d=2,3$, 
but not for $d=4$. 

\begin{Proposition} \sl Assume $d=2$. Then 
\[ J(A,M) \iff (A,M)_2=0. \] 
\end{Proposition}
\demo 
The forward implication is clear. 
For the converse, assume $(A,M)_2=0$. Then 
$\Gordan{A}{M}{A}{2}{2}{2}{0}{1}{1}$ implies that 
$((A,M)_1,A)_1 = -\frac{1}{2} \, (A,A)_2 \, M$. 
If $(A,A)_2 \neq 0$, then let 
\[ B = \frac{2}{(A,A)_2} \, (A,M)_1. \] 

If $(A,A)_2 =0$, then by a change of variable, 
we may assume $A = x_0^2$. Then $(A,M)_2 = 0$ implies 
that $M = c_1 \, x_0^2 + c_2 \, x_0 x_1$. Now 
let $B = c_1 \, x_0 x_1 + \frac{c_2}{2} \, x_1^2$. 
In either case, $(A,B)_1=M$. \qed 

\begin{Proposition} \sl 
Assume $d=3$, then 
\[ J(A,M) \iff ((A,M)_2,A)_1=0. \] 
\end{Proposition} 
\demo By Lemma \ref{lemma.jac}, $((A,M)_2,A)_1=0$ iff $(A,M)_2 = \cons \, A$. 
This shows the forward implication. 

Conversely, assume that $(A,M)_2 = c \, A$ for some 
constant $c$. I claim that the map 
\[ \psi_{M,6 \, c}: S_3 \lra S_3, \quad 
F \lra (M,F)_2 - c \, F \] 
is skew-symmetric. Indeed, 
\[ \delta_{\psi_{M,6 \, c}(F)}(G) 
= ((M,F)_2,G)_3 - c \, (F,G)_3. \] 
Using $\Gordan{M}{F}{G}{4}{3}{3}{1}{2}{2}$, this can be 
transformed into 
\[ -((M,G)_2,F)_3 + c \, (G,F)_3 = - \delta_{\psi_{M,6 \, c}(G)}(F).
\] 
This proves the claim, and implies that the rank of $\psi_{M,6 \, c}$ 
must be even. Suppose that $A$ and another form $B$ span its 
kernel. Then by Proposition \ref{prop.pencil}, $(A,B)_1 = M$ 
(after multiplying $B$ by a constant if necessary). \qed 

\begin{Example} \rm 
Assume $d=4$, and let $A = (x_0x_1)^2, M = (x_0x_1)^3$. Then 
$A$ divides $(A,M)_2 = \cons \,(x_0x_1)^3$. However there exists no 
$B$ such that $(A,B)_1 = M$. Indeed, 
\[ (A,B)_1 = \cons  x_0 \, x_1(x_1 B_{x_1} - x_0 B_{x_0}) =  
(x_0 \, x_1)^3 \] 
would imply $x_1 B_{x_1} - x_0 B_{x_0} = \cons  (x_0x_1)^2$. 
But then $B = \cons (x_0x_1)^2$, which is absurd. 
\end{Example} 

The two propositions above suggest the following natural problem: 
\begin{Problem} \rm 
Find a (finite) number of joint covariants of $A,M$ 
which simultaneously vanish iff $J(A,M)$ holds. 
\end{Problem} 

\subsection{The resultant} 
Let $\Res(A,B)$ denote the resultant of $A,B$. Up to a scalar, it is 
equal to the discriminant of $\R(A,B)$ (regarded as a quadratic form). 
Since the latter implicitly depends only on $T_1,T_3$, the 
following problem is natural: 
\begin{Problem} \rm 
Give an explicit formula (in a reasonable sense) 
for $\Res(A,B)$ as a joint invariant of $T_1$ and $T_3$. 
\end{Problem} 
For instance, if $d=2$ then $\cons \Res(A,B) = (T_1,T_1)_2$. 
\begin{Proposition} \sl If $d=3$, then 
\[ \cons \Res(A,B) = T_3 \, (T_1,T_1)_4 - 6 \, (T_1,(T_1,T_1)_2)_4. \] 
\end{Proposition}
\demo 
By construction, $\Res = \Res(A,B)$ is joint invariant of total degree $3$ 
in $T_1,T_3$. Every joint invariant is a linear combination of 
compound transvectants (see \cite[p.~92]{Glenn}), hence $\Res$ is a 
linear combination of terms of the form 
\[ (X_1,(X_2,X_3)_a)_b, \] 
where $a,b$ are integers, and each $X_i$ stands for 
either $T_1$ or $T_3$. Since the total order must be zero, 
$\sum\limits_i \ord X_i = 2(a+b)$. Using 
properties (\ref{trans.1}),(\ref{trans.2}), we are left with only 
two possibilities, namely 
\[ (T_3,(T_1,T_1)_4)_0, \quad (T_1,(T_1,T_1)_2)_4. \] 
Now specialize to 
$A = x_0 x_1(x_0-x_1), B = x_0(x_0 + x_1)(x_0+2x_1)$ and 
use the method of undetermined coefficients. \qed 

\subsection{The `minimal equation' for $T_3$} 
Consider the following equivalence relation on pairs $(A,B)$ of 
independent order $d$ forms: 
\[ (A,B) \sim (\alpha A + \beta B, \gamma A + \delta B)
\quad \text{if $\alpha \delta - \beta \gamma =1$.} 
\] 
An equivalence class determines and is determined by $T_1,T_3$. 
Let $\F$ denote the set of equivalence classes, and 
consider the map 
\[ \pi: \F \lra {\mathbb A}^{2d-1}, \quad 
(A,B) \lra T_1. \] 
It is known that $\pi$ has finite fibres, and the cardinality of 
the general fibre is equal to the Catalan number 
$\rho(d) = \frac{1}{d}\binom{2d-2}{d-1}$ 
(see \cite[Theorem 1.3]{Goldberg1}). 

Now assume $d=4$, then $\rho(4) =5$. Let $A,B$ be forms of order $4$ 
with {\sl indeterminate} coefficients, and write 
\[ T_1 = \sum_i \binom{6}{i} u_i \, x_0^{6-i}x_1^i, \quad 
   T_3 = \sum_j \binom{2}{j} v_j \, x_0^{2-j}x_1^j, 
\] 
where $u_i,v_j$ are functions of the coefficients of $A,B$. 
The map $\pi$ corresponds to a degree $5$ field extension 
$K \subseteq L$, where 
\[ K = \complex(u_0,\dots,u_6), 
\quad L = K(v_0,v_1,v_2).\] 

We recall the concept of a seminvariant of a form: it is an expression 
in the coefficients of the form which remains unchanged by a substitution 
\begin{equation} 
x_0 \lra x_0 + c \, x_1, \; x_1 \lra x_1; \quad c \in \complex. 
\label{subsU} \end{equation}
An alternative is to define it as the leading coefficient 
of a covariant (see \cite[\S 32]{GrYo}). Let 
\begin{equation} 
v_0^5 + \sum\limits_{i=1}^5 \, l_i \, 
v_0^{5-i} = 0, \qquad l_i \in K,
\label{mineq} \end{equation}
denote the unique minimal equation of $v_0$ over $K$. 
Firstly, since $v_0$ is a seminvariant of $T_3$ 
and substitutions in (\ref{subsU}) must leave (\ref{mineq}) unchanged, 
all the $l_i$ are seminvariants of $T_1$. 
Secondly, by the main theorem of \cite[\S 33]{GrYo}, any algebraic relation 
between the seminvariants translates into a relation between the corresponding 
covariants. That is to say, we must have an identity 
\begin{equation} 
 T_3^5 + \sum\limits_{i=1}^5 \, \Lambda_i \, T_3^{\, 5-i} = 0, 
\label{mineq.cov} \end{equation}
where $\Lambda_i$ are covariants of $T_1$, and (\ref{mineq.cov}) reduces 
to (\ref{mineq}) by the substitution $x_0:=1,x_1:=0$. 
By homogeneity, $\Lambda_i$ must have degree-order $(i,2i)$. 
\subsection{} 
A complete set of generators for the ring of covariants of 
order $6$ forms is given in \cite[\S 134]{GrYo}. It is then  
a routine matter to identify the $\Lambda_i$ by the method 
of undetermined coefficients. I omit all calculations and 
merely state the result. Define the following covariants of $T_1$. 
\[ \begin{array}{lll} 
q_{20} = (T_1,T_1)_6, & q_{24} = (T_1,T_1)_4, & q_{28} = (T_1,T_1)_2, \\ 
q_{32} = (T_1,q_{24})_4, & q_{36} = (T_1,q_{24})_2, & 
q_{38} = (T_1,q_{24})_1, \\ 
q_{44} = (T_1,q_{32})_2. 
\end{array} \] 
These are all taken from the table in \cite[p.~156]{GrYo}, 
but the notation is modified so that $q_{ab}$ is of 
degree-order $(a,b)$. There can be no covariant of degree-order $(1,2)$, 
hence $\Lambda_1 = 0$. The others are 
\[ \begin{aligned} 
\Lambda_2 & = -\frac{125}{8} \, q_{24} \\ 
\Lambda_3 & = \frac{625}{24} \, q_{36} + \frac{125}{36} \, T_1 \, q_{20} \\ 
\Lambda_4 & = \frac{3125}{48} \, q_{24}^2 -\frac{625}{96} \, q_{20} \, q_{28} 
         -\frac{3125}{96} \, T_1 \, q_{32} \\ 
\Lambda_5 & = \frac{3125}{64} \, T_1 \, q_{44} +
         \frac{3125}{64} \, q_{32} \, q_{28} 
        -\frac{3125}{16} \, q_{36} \, q_{24}
        -\frac{3125}{192} \, T_1 \, q_{20} \, q_{24}.
\end{aligned} \] 

\smallskip 

\begin{Problem} \rm 
Find the equation analogous to (\ref{mineq.cov}) for arbitrary $d$. 
\end{Problem}

\bibliographystyle{plain} 

\vspace{1cm} 

\parbox{10cm}{\small 
Jaydeep V.~Chipalkatti \\ 
Department of Mathematics \\ 
University of British Columbia \\ 
Vancouver BC V6T 1Z2, Canada. \\ 
{\tt jaydeep@math.ubc.ca}}
\end{document}